# Quantization of Fractional Singular Lagrangian Systems with Second-Order Derivatives Using Path Integral Method

Eyad Hasan Hasan and Osama Abdalla Abu-Haija

Tafila Technical University, Faculty of Science, Applied Physics Department, P. O. Box: 179, Tafila 66110, Jordan

**Abstract**

We examined the fractional second order singular Lagrangian systems. We wrote the action principal function and equations of motion as fractional total differential equations. Also, we constructed the set of Hamilton–Jacobi partial differential equations (HJPDEs) within fractional calculus. We formulated the fractional path integral quantization for these systems. A mathematical example is examined with first- and second-class constraints.

**Email**: dr_eyad2004@ttu.edu.jo, iyad973@yahoo.com

## 1. Introduction

The quantization of singular Lagrangian systems has been studied with increasing interest and treated first by Dirac [1, 2], for quantizing physical systems. Following Dirac's, researchers have developed a formalism for investigating these systems using the canonical formalism [3]. The set of (HJPDEs) is constructed, and they wrote equations of motion as total differential equations. Then, the WKB approximation and path integral technique are quantized by using this formalism [4-8].

Fractional calculus with singular systems has been treated with more interests and importance [9-18]. Recently, the Euler-Lagrange formalism is analyzed for second order Lagrangian systems within fractional calculus and the fractional Hamilton-Jacobi formalism for these systems is discussed [15, 16]. More recently, authors have constructed a formalism using the canonical method for quantizing singular systems using the path integral approach and WKB approximation for first-order derivatives [17, 18]. In this paper, we would like to extend our work for Lagrangians having second-order derivatives.

The most important definitions of fractional calculus [9] are.

(i) The left Riemann–Liouville fractional derivative

$$ {}_aD_t^{\alpha} f(t) = \frac{1}{\Gamma(n-\alpha)} \left( \frac{d}{dt} \right)^n \int_a^t (t-\tau)^{n-\alpha-1} f(\tau) d\tau \,. \tag{1} $$

(ii) The right Riemann–Liouville fractional derivative

$$ {}_tD_b^{\alpha} f(t) = \frac{1}{\Gamma(n-\alpha)} \left( -\frac{d}{dt} \right)^n \int_t^b (\tau-t)^{n-\alpha-1} f(\tau) d\tau \,. \tag{2} $$

where $\Gamma$ is the Euler gamma function, $n \in N$, $n-1 \leq \alpha < n$, $\alpha$ is an integer and these derivatives can be defined as follows:

$$ {}_aD_t^{\alpha} f(t) = \left( \frac{d}{dt} \right)^{\alpha} f(t) \,, \qquad {}_tD_b^{\alpha} f(t) = \left( -\frac{d}{dt} \right)^{\alpha} f(t) \,, \tag{3} $$

**Definition:** Given a function $f : [0, \infty) \to \Re$. Then for all $t > 0$ , $\alpha \in (0,1)$, let

$$ D^{\alpha}(f)(t) = \lim_{\varepsilon \to 0} \frac{f(t + \varepsilon t^{1-\alpha}) - f(t)}{\varepsilon} \,. \tag{4} $$

$D^{\alpha}$ is called the conformal fractional derivative of $f$ of order of $\alpha$ [17].

In this work, we aim to construct the formalism for quantizing singular Lagrangians systems with second-order derivatives in fractional form.

## 2. Fractional Path Integral Quantization and Fractional Second-order Singular Lagrangian and

Following Hasan [19], we will use a formalism for second-order fractional singular Lagrangian systems to be applicable for quantizing these systems using the path integral approach. The Lagrangian formalism of second-order in fractional derivatives is given by [19]

$$ L = L(D^{\alpha-1} q_i, D^{\alpha} q_i, D^{2\alpha} q_i, t) \,. \tag{5} $$

Where $D^{\alpha} q_i$ are the conformal fractional derivatives of the coordinates $q_i$ [17].

The Lagrangian and Hamiltonian formalism for second-order derivatives have been studied by Ostrogradski [ 20] and these derivatives have been treated as coordinates. Therefore, we can treat these derivatives $D^{\alpha-1}q_i$ and $D^{\alpha}q_i$ as coordinates. Thus, the Poisson brackets can be defined as

$$\{A,B\} \equiv \frac{\partial A}{\partial D^{\alpha-1}q_i}\frac{\partial B}{\partial p_i} - \frac{\partial A}{\partial p_i}\frac{\partial B}{\partial D^{\alpha-1}q_i} + \frac{\partial A}{\partial D^{\alpha}q_i}\frac{\partial B}{\partial \pi_i} - \frac{\partial A}{\partial \pi_i}\frac{\partial B}{\partial D^{\alpha}q_i} \ . \tag{6}$$

Here, the functions A and B are described in term of the canonical variables $D^{\alpha-1}q_i$, $D^{\alpha}q_i$, $p_i$ and $\pi_i$. Thus, the generalized momenta $p_i$ and $\pi_i$ are conjugated to the generalized coordinates $D^{\alpha-1}q_i$ and $D^{\alpha}q_i$ respectively.

Now, the fractional of the Hessian matrix is defined as [19]

$$W_{ij} = \frac{\partial^2 L}{\partial D^{2\alpha}q_i \partial D^{2\alpha}q_j} \qquad i,j = 1,2,\ldots,N \tag{7}$$

The fractional Lagrangian is called regular if it's rank is $N$ otherwise the Lagrangian is singular $N-R, \quad R<N$ . Dirac showed in his formalism for investigating Lagrangians having singular nature that the number of degrees of freedom can be reduced from $N$ to $N-R$ due to the constraints [1, 2]. Thus, we can define the momenta $\pi_i$ conjugated to the coordinates $D^{\alpha}q_i$ as [19]:

$$\pi_a = \frac{\partial L}{\partial D^{2\alpha}q_a}, \qquad a = 1.2,\ldots,N-R \tag{8}$$

$$\pi_\mu = \frac{\partial L}{\partial D^{2\alpha}q_\mu} \ . \qquad \mu = N-R+1,\ldots,N \ . \tag{9}$$

Also, the momenta $p_i$ conjugated to the coordinates $D^{\alpha-1}q_i$ can be defined as [19]:

$$p_a = \frac{\partial L}{\partial D^{\alpha}q_a} - \frac{d}{dt}\left(\frac{\partial L}{\partial D^{2\alpha}q_a}\right); \tag{10}$$

$$p_{\mu} = \frac{\partial L}{\partial D^{\alpha} q_{\mu}} - \frac{d}{dt}\left( \frac{\partial L}{\partial D^{2\alpha} q_{\mu}} \right). \tag{11}$$

Where

$$\pi_{\mu} = -H_{\mu}^{\pi}(D^{\alpha-1} q_i, D^{\alpha} q_i, p_a, \pi_a) \tag{12}$$

and

$$p_{\mu} = -H_{\mu}^{p}(D^{\alpha-1} q_i, D^{\alpha} q_i, p_a, \pi_a) \tag{13}$$

Thus, Equations (12) and (13) represent primary constraints [1, 2] and can be written as

$$H'^{p}_{\mu}(D^{\alpha-1} q_i, D^{\alpha} q_i, p_i, \pi_i) = p_{\mu} + H_{\mu}^{p} = 0 \quad ; \tag{14}$$

$$H'^{\pi}_{\mu}(D^{\alpha-1} q_i, D^{\alpha} q_i, p_i, \pi_i) = \pi_{\mu} + H_{\mu}^{\pi} = 0 \quad . \tag{15}$$

We can calculate the Hamiltonian $H_{\circ}$ in fractional form as

$$H_{\circ} = -L(D^{\alpha-1} q_i, D^{2\alpha} q_{\mu}, D^{\alpha} q_i, D^{2\alpha} q_a) + p_a D^{\alpha} q_a + \pi_a D^{2\alpha} q_a - D^{\alpha} q_{\mu} H_{\mu}^{p} - D^{2\alpha} q_{\mu} H_{\mu}^{\pi}$$

$$\mu = 1, \ldots\ldots, R \; ; \qquad a = R+1, \ldots, N \, . \tag{16}$$

A natural of singular Lagrangian indicates that the generalized momenta $p_{\mu}$ and $\pi_{\mu}$ are not independent of $p_a$ and $\pi_a$. Thus, we can write the set of Hamilton-Jacobi partial differential equations as

$$H'_{\circ} = p_{\circ} + H_{\circ}\left(D^{\alpha-1}q_i, D^{\alpha}q_i, \frac{\partial S}{\partial D^{\alpha-1}q_a}, \frac{\partial S}{\partial D^{\alpha}q_a}\right) = 0\,; \tag{17a}$$

$$H'^{p}_{\mu} = p_{\mu} + H^{p}_{\mu}\left(D^{\alpha-1}q_i, D^{\alpha}q_i, \frac{\partial S}{\partial D^{\alpha-1}q_a}, \frac{\partial S}{\partial D^{\alpha}q_a}\right) = 0\,; \tag{17b}$$

$$H'^{\pi}_{\mu} = \pi_{\mu} + H^{\pi}_{\mu}\left(D^{\alpha-1}q_i, D^{\alpha}q_i, \frac{\partial S}{\partial D^{\alpha-1}q_a}, \frac{\partial S}{\partial D^{\alpha}q_a}\right) = 0. \tag{17c}$$

Here, $S = S(D^{\alpha-1}q_a, D^{\alpha-1}q_{\mu}, D^{\alpha}q_a, D^{\alpha}q_{\mu}, t)$ is the fractional Hamilton's principal function.

Considering the definitions of the generalized momenta as

$p_a = \dfrac{\partial S}{\partial D^{\alpha-1}q_a}$, $p_{\mu} = \dfrac{\partial S}{\partial D^{\alpha-1}q_{\mu}}$ , $\pi_a = \dfrac{\partial S}{\partial D^{\alpha}q_a}$, $\pi_{\mu} = \dfrac{\partial S}{\partial D^{\alpha}q_{\mu}}$ and $p_{\circ} = \dfrac{\partial S}{\partial t}$.

Researchers in reference [5] wrote the Hamilton's principal function and the equations of motion as total differential equations, also we can write these equations in fractional form as follows [5]:

$$dD^{\alpha-1}q_a = \frac{\partial H'_{\circ}}{\partial p_a}dt + \frac{\partial H'^{p}_{\mu}}{\partial p_a}dD^{\alpha-1}q_{\mu} + \frac{\partial H'^{\pi}_{\mu}}{\partial p_a}dD^{\alpha}q_{\mu}, \tag{18}$$

$$dD^{\alpha}q_a = \frac{\partial H'_{\circ}}{\partial \pi_a}dt + \frac{\partial H'^{p}_{\mu}}{\partial \pi_a}dD^{\alpha-1}q_{\mu} + \frac{\partial H'^{\pi}_{\mu}}{\partial \pi_a}dD^{\alpha}q_{\mu}, \tag{19}$$

$$-dp_i = \frac{\partial H'_{\circ}}{\partial D^{\alpha-1}q_i}dt + \frac{\partial H'^{p}_{\mu}}{\partial D^{\alpha-1}q_i}dD^{\alpha-1}q_{\mu} + \frac{\partial H'^{\pi}_{\mu}}{\partial D^{\alpha-1}q_i}dD^{\alpha}q_{\mu}, \tag{20}$$

$$-d\pi_i = \frac{\partial H'_{\circ}}{\partial D^{\alpha}q_i}dt + \frac{\partial H'^{p}_{\mu}}{\partial D^{\alpha}q_i}dD^{\alpha-1}q_{\mu} + \frac{\partial H'^{\pi}_{\mu}}{\partial D^{\alpha}q_i}dD^{\alpha}q_{\mu}, \tag{21}$$

$$dS = \left(-H_{\circ} + p_a\frac{\partial H'_{\circ}}{\partial p_a} + \pi_a\frac{\partial H'_{\circ}}{\partial \pi_a}\right)dt + \left(-H^{p}_{\mu} + p_a\frac{\partial H'^{p}_{\mu}}{\partial p_a} + \pi_a\frac{\partial H'^{p}_{\mu}}{\partial \pi_a}\right)dD^{\alpha-1}q_{\mu} + \left(-H^{\pi}_{\mu} + p_a\frac{\partial H'^{\pi}_{\mu}}{\partial p_a} + \pi_a\frac{\partial H'^{\pi}_{\mu}}{\partial \pi_a}\right)dD^{\alpha}q_{\mu}. \tag{22}$$

The set of equations (18-22) are integrable if the total derivative of equation (17) is zero [3, 5]

$$dH'_{\circ} = 0\,;\ dH'^{p}_{\mu} = 0\,;\ dH'^{\pi}_{\mu} = 0. \tag{23}$$

Thus, the degrees of freedom are reduced from $N$ to $N-R$, and the canonical phase space coordinates have been reduced from $\{D^{\alpha-1}q_i, p_i, D^{\alpha}q_i, \pi_i\}$ to $\{D^{\alpha-1}q_a, p_a, D^{\alpha}q_a, \pi_a\}$. Therefore, the path integral approach can be represented in the fractional form as

$$K(D^{\alpha-1}q_a, D^{\alpha-1}q_\mu, D^{\alpha}q_a, D^{\alpha}q_\mu, t) = \int \prod_{a=1}^{N-R} dD^{\alpha-1}q_a dD^{\alpha}q_a dp_a d\pi_a \exp i \begin{bmatrix} \int\left(-H_{\circ} + p_a \dfrac{\partial H'_{\circ}}{\partial p_a} + \pi_a \dfrac{\partial H'_{\circ}}{\partial \pi_a}\right)dt + \\ \int \left(-H^{p}_{\mu} + p_a \dfrac{\partial H'^{p}_{\mu}}{\partial p_a} + \pi_a \dfrac{\partial H'^{p}_{\mu}}{\partial \pi_a}\right) dD^{\alpha-1}q_\mu + \\ \int\left(-H^{\pi}_{\mu} + p_a \dfrac{\partial H'^{\pi}_{\mu}}{\partial p_a} + \pi_a \dfrac{\partial H'^{\pi}_{\mu}}{\partial \pi_a}\right) dD^{\alpha}q_\mu \end{bmatrix}$$

. (24)

## 4. Example:

We will discuss an example of second-order fractional singular Lagrangian has primary and secondary constraints:

$$L = \frac{1}{2}\left((D^{2\alpha}q_1)^2 + (D^{2\alpha}q_2)^2\right) - \frac{1}{2}[(D^{\alpha}q_1)^2 + (D^{\alpha}q_2)^2] + \frac{1}{2}D^{\alpha-1}q_3 + D^{\alpha}q_3 D^{2\alpha}q_3. \quad (25)$$

The generalized momenta read.

$$p_1 = -D^{\alpha}q_1 - D^{3\alpha}q_1 \quad ; \quad (26a)$$

$$p_2 = -D^{\alpha}q_2 - D^{3\alpha}q_2; \quad (26b)$$

$$p_3 = 0 = -H^{p}_{3}; \quad (26c)$$

$$\pi_1 = D^{2\alpha}q_1; \quad (26d)$$

$$\pi_2 = D^{2\alpha}q_2; \quad (26e)$$

$$\pi_3 = D^{\alpha}q_3 = -H^{\pi}_{3}. \quad (26f)$$

Here, Equations (26c) and (26f) can be written as

$$H_3'^{p} = p_3 = 0\,; \qquad (27)$$

$$H_3'^{\pi} = \pi_3 - D^{\alpha} q_3 = 0\,. \qquad (28)$$

and represent as primary constraints [1,2].

The Hamiltonian $H_0$ is calculated as

$$H_{\circ} = p_1 D^{\alpha} q_1 + p_2 D^{\alpha} q_2 + \frac{1}{2}[(D^{\alpha} q_1)^2 + (D^{\alpha} q_2)^2] - \frac{1}{2}(D^{\alpha-1} q_3)^2 + \frac{1}{2}(\pi_1^2 + \pi_2^2). \qquad (29)$$

The set of HJPDEs, can be written as.

$$H_{\circ}' = p_{\circ} + H_{\circ} = p_{\circ} + p_1 D^{\alpha} q_1 + p_2 D^{\alpha} q_2 + \frac{1}{2}[(D^{\alpha} q_1)^2 + (D^{\alpha} q_2)^2] - \frac{1}{2}(D^{\alpha-1} q_3)^2 + \frac{1}{2}(\pi_1^2 + \pi_2^2). \qquad (30)$$

$$H_3'^{p} = p_3 = 0\,; \qquad (31)$$

$$H_3'^{\pi} = \pi_3 - D^{\alpha} q_3 = 0\,. \qquad (32)$$

Using the fundamental Poisson brackets

$\{D^{\alpha-1} q_i, p_j\} \equiv \delta_{ij}$ and $\{D^{\alpha} q_i, \pi_j\} \equiv \delta_{ij}$,

$\{D^{\alpha-1} q_i, D^{\alpha-1} q_j\} = \{D^{\alpha} q_i, D^{\alpha} q_j\} = 0 = \{D^{\alpha} q_i, D^{\alpha-1} q_j\} = \{p_i, \pi_i\}$, where

$i, j = 1, \ldots, N$

Thus, the Poisson bracket of $H_3'^{p}$ and $H_0'$

$\{H_3'^{p}, H_{\circ}'\} = D^{\alpha-1} q_3 = H_3''^{p}$ is not identically zero.

it gives a new (secondary) constraint [1, 2]:

$$H_3''^{p} = D^{\alpha-1} q_3 = 0\,. \qquad (33)$$

Following Dirac's classification [1, 2], the constraint (32) is first- and (31) and (33) are second-class. There are no further constraints.

The equations of motion (18- 22) can be calculated as

$$dD^{\alpha-1}q_1 = D^{\alpha}q_1 dt, \tag{34}$$

$$dD^{\alpha-1}q_2 = D^{\alpha}q_2 dt, \tag{35}$$

$$dD^{\alpha}q_1 = \pi_1 dt, \tag{36}$$

$$dD^{\alpha}q_2 = \pi_2 dt, \tag{37}$$

$$dp_1 = 0, \tag{38}$$

$$dp_2 = 0, \tag{39}$$

$$dp_3 = 0, \tag{40}$$

$$-d\pi_1 = (p_1 + D^{\alpha}q_1)dt, \tag{41}$$

$$-d\pi_2 = (p_2 + D^{\alpha}q_2)dt, \tag{42}$$

$$d\pi_3 = dD^{\alpha}q_3, \tag{43}$$

$$dS = \frac{1}{2}([\pi_1^2 + \frac{1}{2}\pi_2^2] - [(D^{\alpha}q_1)^2 + (D^{\alpha}q_2)^2])dt + \pi_3 dD^{\alpha}q_3 \tag{43}$$

Considering $\pi_3 = D^{\alpha}q_3$

$$dS = \frac{1}{2}([\pi_1^2 + \frac{1}{2}\pi_2^2] - [(D^{\alpha}q_1)^2 + (D^{\alpha}q_2)^2])dt + D^{\alpha}q_3 dD^{\alpha}q_3 . \tag{44}$$

Integrate Eq. (44), the action function becomes.

$$S = \int \frac{1}{2}([\pi_1^2 + \frac{1}{2}\pi_2^2] - [(D^{\alpha}q_1)^2 + (D^{\alpha}q_2)^2])dt + \frac{1}{2}(D^{\alpha}q_3)^2 . \tag{45}$$

Finally, by obtaining the fractional action function $S$, we can represent the path integral approach in fractional form as

$$\begin{aligned} K(D^{\alpha-1}q_1, D^{\alpha-1}q_2, D^{\alpha-1}q_3, D^{\alpha}q_1, D^{\alpha}q_2, D^{\alpha}q_3, t) &= \int dD^{\alpha-1}q_1 dD^{\alpha-1}q_2 dD^{\alpha}q_1 dD^{\alpha}q_2 dp_1 dp_2 d\pi_1 d\pi_2 \\ \exp i\left[\int \frac{1}{2}([\pi_1^2 + \pi_2^2] - [(D^{\alpha}q_1)^2 + (D^{\alpha}q_2)^2])dt + \frac{1}{2}\pi_3^2\right] \end{aligned} . \tag{46}$$

Taking into account, $D^{\alpha-1}q_3 = 0$ , $\pi_3 = D^{\alpha}q_3$ thus, $D^{\alpha}q_3 = 0$
Thus, equation (46) becomes.

$$K(D^{\alpha-1}q_1, D^{\alpha-1}q_2, D^{\alpha-1}q_3, D^{\alpha}q_1, D^{\alpha}q_2, D^{\alpha}q_3, t) = \int dD^{\alpha-1}q_1 dD^{\alpha-1}q_2 dD^{\alpha}q_1 dD^{\alpha}q_2 dp_1 dp_2 d\pi_1 d\pi_2$$
$$\exp i\left[\frac{1}{2}\int ([\pi_1^2 + \pi_2^2] - [(D^{\alpha}q_1)^2 + (D^{\alpha}q_2)^2])dt\right] . \qquad (47)$$

## Conclusion

In this work, we constructed a formalism for quantizing singular systems using fractional path integral technique. We wrote Hamilton's principal function and equations of motion in fractional form as total differential equations. Calculating the fractional Hamilton's principal function enabled us to formulate the fractional path integral technique. Then, the quantization can be carried out. Finally, we discussed a mathematical example to demonstrate our formalism.